\documentclass{amsart}

\usepackage{amsmath,amsthm,amsbsy,amssymb}
\usepackage{enumitem}
\usepackage{xcolor} 
\definecolor{orange}{rgb}{1,0.5,0}
\definecolor{Red}{rgb}{.795,0.015,0.017}
\definecolor{Ggreen}{rgb}{0.,0.675,0.0128}
\definecolor{Bblue}{rgb}{0.16,.32,0.91}
 \usepackage[backref=page]{hyperref}
\hypersetup{
    colorlinks = true,
linkcolor={red},
urlcolor={blue},
citecolor={Ggreen},    
urlcolor = {blue},
citebordercolor = {0.33 .58 0.33},
 linkbordercolor = {0.99 .28 0.23},
 pdftitle={A statistical model for points expanding in higher dimensions while being tied to bijective involutions},
 pdfauthor={CNZ}
}
\usepackage{mathabx}
\usepackage[mathscr]{euscript}
\newcommand{\PI}{\mathscr{P}(\cI)}
\newcommand{\PIp}{\mathscr{P}(\cI')}
\newcommand{\PIIp}{\mathscr{P}(\cI\cup\cI')}

\newcommand{\Fpsi}{\mathscr{F}_\psi}

\usepackage{mathrsfs}% scriptfont

\newcommand{\cB}{\mathcal B}

\newcommand{\cE}{\mathcal E}
\newcommand{\cF}{\mathcal F}

\newcommand{\cI}{\mathcal I}

\newcommand{\cM}{\mathcal M}

\newcommand{\cP}{\mathcal P}

\newcommand{\cS}{\mathcal S}
\newcommand{\cT}{\mathcal T}

\newcommand{\cX}{\mathcal X}

\usepackage{bm}
\newcommand*{\B}[1]{\ifmmode\bm{#1}\else\textbf{#1}\fi}

\newcommand{\bx}{\B{x}}

\newcommand{\bXM}{{\B{\cX}}_{\!\psi}}

\def\ZZ{\mathbb{Z}}

\newcommand{\sdfrac}[2]{\mbox{\small$\displaystyle\frac{#1}{#2}$}}
\newcommand{\fdfrac}[2]{\mbox{\footnotesize$\displaystyle\frac{#1}{#2}$}}

\DeclareMathOperator{\A}{A}
\DeclareMathOperator{\M}{M_2}

\usepackage{float}
\renewcommand{\pmod}[1]{\left( \mathrm{ mod\;}#1\right)}
\newtheorem{theorem}{Theorem}
\newtheorem{corollary}{Corollary}

\newtheorem{lemma}{Lemma}[section]

\theoremstyle{remark}
\newtheorem{remark}{Remark}[section]
\newtheorem*{remark*}{Remark}

\theoremstyle{definition}

\renewcommand*{\backref}[1]{}
\renewcommand*{\backrefalt}[4]{%
  \ifcase #1 %
No citations.% use \relax if you do not want the "No citations" message
  \or
(page #2).%
  \else
(pages #2).%
  \fi%
}
\usepackage{cite}

\newtheorem*{oldproof}{Proof}

\title[A statistical model for points---tied to bijective involutions]{A statistical model for points expanding in higher dimensions while being tied to bijective involutions}

\author{Cristian Cobeli}
\address{CC: ``Simion Stoilow''\\ 
Institute of Mathematics\\
of the Romanian Academy\\
21 Calea Griviței Street \\
P. O. Box 1-764 \\
Bucharest 014700 \\
Romania\\
cristian.cobeli@imar.ro\\}

\author{The Nguyen}
\address{TN: Department of Mathematics\\
University of Illinois at Urbana-Champaign\\
Altgeld Hall, 1409 W. Green Street\\
Urbana, IL, 61801, USA\\
thevn2@illinois.edu}

\author{Alexandru Zaharescu}
\address{AZ: Department of Mathematics\\
University of Illinois at Urbana-Champaign\\
Altgeld Hall, 1409 W. Green Street\\
Urbana, IL, 61801, USA\\
and\\
``Simion Stoilow''\\ Institute of Mathematics of the Romanian Academy\\
21 Calea Griviței Street \\
P. O. Box 1-764, Bucharest 014700 \\
Romania\\
zaharesc@illinois.edu}

\makeatletter
\@namedef{subjclassname@2020}{\textup{2020} Mathematics Subject Classification}
\makeatother
\subjclass[2020]{primary 11N69; secondary 49J55, 65C10}
\keywords{Statistical model; modular mappings; 
random sequences; Poisson distribution; fixed points; Stauduhar's Conjecture}

\begin{document}

\begin{abstract}
Let $\cM$ be a set with $M$ elements, 
let $\psi :\cM\to\cM$ be a bijective involution,
and let~$\bXM$ be the set of sequences $(x_1,\dots,x_M)\in\cM^M$
with the property that \mbox{$x_{M+1-j} = \psi(x_j)$} for $1\le j\le M$.
This framework can be used to infer the possible
distribution of sequences, such as the modular ones, 
that pose challenges for conventional methods.

We prove that when $M$ is even, there exists a limit probability density function
that weighs the parameter $k$ that counts the appearances of the  
elements of $\cM$ among the terms of sequences $\bx\in\bXM$. 
It turns out that the number of fixed points of $\psi$ influences the 
probability density function, which decomposes into two
pieces, each multiplied by complementary factors, and the smaller of the two pieces appears only when $k$ is even. 

Applying the model, we find a threshold from which almost all sequences contain related terms with prescribed frequencies.
\end{abstract}

\maketitle

\section{Introduction}

%%%%%%%%%%%%%%%%%%%%%%%%%%%%%%%%%
We introduce and study a statistical model aimed primarily at understanding the expected distribution of sequences of numbers which, despite being a common occurrence, the actual proof of the observed properties is beyond the possibilities offered by known algebraic, analytical, and probabilistic techniques.

An illustrative case is that of sequences composed of residue classes $Q(n)$ modulo a prime number $p$. 
In this case, $Q(n)$ can be a power, a polynomial, or a rational function.
Or, in a more intriguing example resembling randomness, the `mixing diagonal product' of the factorials 
\mbox{$1!, 2!, \dots,(p-1)!\pmod p$}.
Based on standard expectations (see
Kolchin et al.~\cite{KSC1978},
Mikhajlov~\cite{Mik1980},
Vatutin and Mikhajlov~\cite{VM1982}), 
it is conjectured that this sequence follows 
the distribution of a Poisson process 
(see Stauduhar~\cite[Problem 7]{NTC1963}, 
Guy~\cite[Problem F11]{Guy2004} and \cite{CVZ2000, CPZ2024}).
This means that for any integer $k\ge 0$, 
the proportion of residue classes $y\in\ZZ/p\ZZ$
for which there are exactly $k$ integers 
$1 \le x \le p - 1$ for which $x! \equiv y \pmod p$ 
is expected to tend to $1/(ek!)$, as $p$ tends to infinity.
Then, when $k=0$, it would follow that about $1/e$ residue classes are missed by the sequence of factorials $\pmod p$.

But while this is the general behavior of almost all sequences,
without having any specific property to distinguish them, 
obtaining results towards proving that 
the factorials modulo $p$ 
follow a Poisson distribution
is a notable challenge.
Over time, numerous studies have been dedicated to this problem (see Banks et al.~\cite{BLSS2005}, 
Chen and Dai~\cite{CD2006},
Garaev and García~\cite{GG2007},
Garaev and Hern\'andez~\cite{GH2017},
Klurman and Munsch~\cite{KM2017}).
Additionally, other versions of the problem involving products
of factorials (see 
Garaev et al.~\cite{GLS2004},
García~\cite{Gar2007, Gar2008},
Luca and Stănică~\cite{LS2003})
or generalizations in other spaces have been tinkered
(see Garaev et al.~\cite{GLS2005},
Bhargava~\cite{Bha2000}, Lagarias and Yangjit~\cite{LY2024} and the references therein). 
The best result so far is due to 
Grebennikov et al.~\cite{GSSV2024}, 
who showed the certain existence of merely 
$\left(\sqrt{2}+o(1)\right)\sqrt{p}$ distinct factorials $\pmod p$ as $p$ tends to infinity.

Let $M$ be a positive integer, and let $\cM$ be a set of  cardinality $|\cM|=M$.
The special feature of our model is that from the entire set of all sequences (or vectors) of the cartesian product $\cM^M$, 
in the statistical accounting, we will only consider  
a small portion of those sequences that are related 
to the arithmetic structure 
of the residue classes in $\ZZ/p\ZZ$.
And a characterization that provides a fine yet sufficiently generous encompassing filter is the linking of residue classes through their intrinsic symmetries.
One example of such symmetry is the one that associates to $x$ 
its inverse $x^{-1}\pmod p$. 
This operation has an immediate effect on the decomposition
of the set of residue classes 
$\{2, 3,\dots,p-2\}=\cI\cup\cI'$ into two subsets~$\cI$ and $\cI'$, 
one containing half and the other containing the inverses from the first half. The consequence is the well-known characterization 
of prime numbers through Wilson's theorem:
\begin{align*}
   (p-1)! \equiv 1\cdot\prod_{x\in\cI}x\prod_{x\in\cI'}x^{-1}\cdot (p-1)
   \equiv \prod_{x\in\cI}x x^{-1} \cdot (-1)\equiv -1\pmod p.
\end{align*}

%%%%%%%%%%%%%%%%%%%%%%%%%%%%%%%%%%%%
\subsection{Assumptions of the model}
Therefore, the configuration adopted for our model consists of the triplet $(\cM,\psi, \bXM)$, where~$\cM$ is a set of $M>0$ elements,
$\psi=\psi_{\cM}$ is an involutive permutation $\psi : \cM \to \cM$,
that is, it satisfies condition  $\psi \circ \psi = \mathrm{Id}$, 
and
\begin{equation}\label{eqDefinitionbXM}
  \bXM := \big\{\bx = (x_1,\dots, x_M) : x_{M+1 - j} = \psi(x_j) \text{ for all } j = 1, 2, \dots, M\big\}.
\end{equation} 

We will assume throughout the following that $M$ is even 
(such as when \mbox{$M=p-1$}, with $p$ prime), 
in order to maintain a straightforward  presentation.
Let us note then that the pool of sequences for our statistics is reduced from all~$M^M$ possibilities to just $|\bXM| = M^{M/2}$.

It turns out that the number of fixed points that $\psi$ could have is proving to play a key factor in our model, 
which is why we are led to assume that there exist 
$\theta\in [0,1]$ and $\eta \in [0,1)$
such that the number of fixed points 
of~$\psi$~is 
% $\theta M\left(1 + O(M^{\eta-1})\right)$ as $M$ tends to infinity.
\begin{align}\label{eqNumberOfFixedPoints}
   |\cF_\psi| = \theta M\left(1 + O(M^{\eta-1})\right),
   \ \ \text{as $M\to\infty$,}
\end{align}
where $\cF_\psi$ is the set of fixed points of $\psi$.

We further let the value of~$\theta$ remain fixed or change with $M$,
for instance, we might have $\theta = a/M$, 
which includes the possibility that $\psi$ could have
a certain number 
$a$ of fixed points, although not necessarily zero.
% We say that a vector $\bx = (x_1, \dots, x_M) \in \bXM$ \textit{misses} the element~$y \in \cM$ if no component of $\bx$ coincides with $y$.

%%%%%%%%%%%%%%%%%%%%%%%%%%%%%%%%%%%%
\subsection{The main results}

The following theorem describes the behavior of the presence of elements $y\in\cM$ in the sequences $\bx\in\bXM$. 
Notably, the values of the probability density function change
based on the parity of the index $k$ of the number of occurrences.
%%%%%%%%%%%%%%%%%%%%%%%%%%%%%%%%%%%%%%%
\begin{theorem}\label{Theorem1}
Let $\cM$ be a set of  cardinality $M$, where $M$ is a positive even integer. Let $\psi=\psi_{\cM}$ be an involutive permutation of $\cM$ and consider the set  $\bXM$ of $\psi$-symmetric vectors of length $M$ defined by~\eqref{eqDefinitionbXM}.
Suppose the number of fixed points of $\psi$
is $\theta M\left(1 + O(M^{\eta-1})\right)$, where
$\theta\in [0,1]$ and $\eta\in [0,1)$.
Let $k$ be a non-negative integer. Then,  when $M$ tends to infinity, we have: 
 \begin{itemize}
     \item[(a)] 
  If $k$ is odd, in asymptotically all~\footnote{Here, ``asymptotically all'' 
  means that the limit of the ratio between the cardinality of 
  the exceptional set of those $\bx\in\bXM$ that do not satisfy 
  the theorem's conclusion and the cardinality of 
  $\bXM$ tends to zero, as $M$ tends to infinity.} 
      vectors $\bx\in\bXM$, there are about
$\frac{(1-\theta)M}{e\,k!}$ distinct components that each occur exactly $k$ times among 
the components of~$\bx$.
     \item[(b)] 
      If $k$ is even, in asymptotically all vectors $\bx\in\bXM$, there are about
$\frac{(1-\theta)M}{e\,k!} + \frac{\theta M}{e^{1/2}\, 2^{k/2} (k/2)!  }$ distinct components that each occur exactly $k$ times among the components of $\bx$.      
 \end{itemize}
\end{theorem}

When $k = 0$, Theorem~\ref{Theorem1} reduces to the following simpler form. 
%%%%%%%%%%%%%%%%%%%%%%%%%%%%%%%%%%%%%%
\begin{corollary}\label{Corollary1}
Let us assume that the hypotheses of 
Theorem~\ref{Theorem1} are satisfied.
Then, when $M \to \infty$, in asymptotically all the 
vectors $\bx\in\bXM$ there are not represented about $\left(\frac{1-\theta}{e} 
 + \frac{\theta}{e^{1/2}}\right) M$  elements of $\cM$. 
\end{corollary}
A precise quantification of the expression in Theorem~\ref{Theorem1} and Corollary~\ref{Corollary1} is given at the end of their proof, in Theorem~\ref{Theorem11} 
and Corollary~\ref{Corollary11}.

Next, we apply the result from Theorem~\ref{Theorem11},
showing that 
if $\psi$ has no fixed points, then
almost all sequences in $\bXM$ have components 
belonging to certain sufficiently large subsets of $\cM$, 
the number of occurrences of these components can be prescribed a priori, and moreover, they are also interlinked by priori fixed bijections.
%%%%%%%%%%%%%%%%%%%%%%%%%%%%%%%%%%%%%5
\begin{theorem}\label{Theorem2}
Let $r, s \ge 0$ be integers such that $r+s \ge 1$.
Then, for any $\varepsilon > 0$, there exists a constant 
$c = c(r,s,\varepsilon) \in (0,1]$ 
and an integer $M_0= M_0(r,s,\varepsilon)>0$ 
such that for any set~$\cM$ of even cardinality $M \ge M_0$  
and any involutive permutation $\psi : \cM \to \cM$ with no fixed points, 
there is a subset 
$\cE = \cE(\cX_{\psi,\cM}) \subset \cX_{\psi,\cM}$ 
with $|\cE| < \varepsilon |\cX_{\psi, \cM}|$ 
such that for any $\mathcal{B} \subset \cM$ 
of size $|\cB| > c |\cM|$,
the following statement holds true.

In any $\bx = (x_1, \dots, x_M) \in \bXM \setminus \cE$
and for any bijection $\phi : \cM \to \cM$, 
there are components $u, v \in \cB\cap\{x_1,\dots,x_M\}$ 
such that: 
    \begin{enumerate}
        \item $u$ appears at most $r$ times in $\bx$; 
        \item $v$ appears at most $s$ times in $\bx$;
        \item %and 
        $\phi(u) = v$.
    \end{enumerate}
\end{theorem}

\begin{theorem}\label{Theorem3}
For all $\varepsilon > 0$, there exists an integer 
$M_0 = M_0(\varepsilon)$ such that for any set~$\cM$ with
even cardinality $|\cM| = M  \ge M_0$ and any 
involutive permutation
% bijection 
$\psi : \cM \to \cM$ with no fixed points
% with $\psi \circ \psi = \mathrm{Id}$, 
there are at least $(1-\varepsilon)|\cX_{\psi, \cM}|$ vectors 
$\bx = (x_1, \dots, x_M)\in\cX_{\psi,\cM}$ 
with the property that
in any subset $\cB \subset \cM$ that is large enough, such that it consists of at least $87\%$ of the elements of~$\cM$,
and for any bijection $\phi : \cM \to \cM$, 
there exists $y = y(\bx, \cB)$ such that 
both~$y$ and $\phi(y)$ are among the components of $\bx$.
\end{theorem}
We note that the percentage $87\%$ in the statement of 
Theorem~\ref{Theorem3} is tight. Indeed, Remark~\ref{Remark86Tight}
at the end of proof of Theorem~\ref{Theorem3} shows that the 
$86\%$ percentage is not sufficient for the statement to remain valid. 
Actually, the precise value of the threshold that makes the statement true or false is $1/2+e^{-1}\approx 0.867879$.
% A = numerical_approx((1/2 + exp(-1)))  #0.867879441171442

%%%%%%%%%%%%%%%%%%%%%%%%%%%%%%%%%
\section{Notation}
Let $\cM$ be a set with cardinality $M$, where  $M\ge 2$ is assumed in all that follows to be even.
For any pair of elements $y,z \in \cM$, we denote 
\begin{align*}
  \delta(y,z) := 
   \begin{cases}
     1, & \text{if } y =z; \\ 
    0, & \text{otherwise.}
   \end{cases}
\end{align*}
Similarly, for any integer $k\ge 0$, if $\bx = (x_1,x_2, \dots, x_M) \in \bXM$ and any $y \in \cM$, we let
\begin{align*}
  \delta_k(\bx,y) := 
   \begin{cases}
     1, & \text{if exactly $k$ components of $\bx$ coincide with $y$}; \\ 
    0, & \text{otherwise.}
   \end{cases}
\end{align*}
Additionally, for any subset $\cI \subset \{1, 2,\dots, M\}$, we denote
\begin{align}\label{eqPixyI}
    \delta_{\cI}(\bx, y) := 
    \begin{cases}
        1 , & \text{ if $x_j = y$ for all $j \in \cI$ and $x_j \neq y$ for all $j \notin \cI$};\\
        0, & \text{otherwise}.
    \end{cases}
\end{align}

For $\bx\in\bXM$ and integers $k=0,1,2,\dots$, 
we keep track of the occurrences of various components in $\bx$ with the following $k$-counter function:
\begin{align*}
  m_k(\bx) := \big|\{y \in \cM : y \text{ is represented exactly $k$ times in } \bx\}\big|.   
\end{align*}

The set of \textit{fixed points} of the bijection $\psi$, which was postulated in the hypothesis from the beginning, plays a central role in the induced `symmetry'.
That is why we let $\Fpsi$ denote the \textit{set of fixed points} of $\psi$, that is, 
\begin{align*}
    \Fpsi  := \{ y \in \cM \,:\, \psi(y) = y\}\,.
\end{align*}
To achieve the outcomes of the model, we will proceed under the assumption that
 there exist $\theta\in [0,1]$ and $\eta \in [0,1)$
such that 
\begin{align*}%\label{eqNumberOfFixedPoints}
   |\cF_\psi| = \theta M\left(1 + O(M^{\eta-1})\right),
\end{align*}
as $M$ tends to infinity.

We use the short notation $[M]$ for the set of integers 
$\{1, 2,\dots, M\}$.

Next, we consider the set of indices of the left leg of the pairs mirrored by $\psi$, where both indices are in a 
subset $\cI\subset \cM$. Thus, we denote: 
\begin{align*}
   \PI := \big\{ j \in [M/2] : 
   \{j, M+1 - j \} \subset \cI\big\}. 
\end{align*}

%%%%%%%%%%%%%%%%%%%%%%%%%%%%%%%%%%%%%%%%%%
\section{\texorpdfstring{The average size of the ratios 
$m_k(\bx)/M$}{The average size of the ratios 
mk(x) over M}}\label{SectionAverage}
We let $\A(k,M)$ denote \textit{the average} value of
the ratios $m_k(\bx)/M$ on all $M$-tuples in $\bXM$:
\begin{align}\label{eqDefAverage}
   \A(k,M) := \frac{1}{|\bXM|}  \sum_{ \bx \in \bXM} \frac{m_k(\bx)}{M}\,. 
\end{align}

By using the $\delta$-functions above, 
the $k$-counter function can be expressed as follows:
\begin{equation}\label{eq:m_k(x)}
  \begin{split}
    m_k(\bx) &= \big|\left\{ y\in \cM \,:\, \delta_k(\bx, y) = 1\right\}\big| \\ 
    &= \sum_{y \in \cM} \delta_k(\bx, y) \\
    &=  \sum_{y \in \cM} 
    \sum_{\substack{\cI \subset [M]\\ |\cI| = k}} 
    \delta_{\cI}(\bx, y)  \\
    &= \sum_{y \in \cM} \sum_{\substack{\cI \subset [M]\\ |\cI| = k}} \prod_{j \in \cI} \delta(x_j,y) \prod_{\ell \in [M]\setminus \cI} \big(1 - \delta(x_\ell, y)\big). 
  \end{split}
\end{equation}

Since the cardinality of $\bXM$ is $|\bXM|=M^{M/2}$, using~\eqref{eq:m_k(x)}, $\A(k,M)$ can be written  as: 
\begin{equation*}% \label{eq:\A(k,M)}
  \begin{split}
  \A(k,M) 
  &= \frac{1}{M^{M/2 + 1}} \sum_{x \in \bXM} \sum_{y \in \cM} \delta_k(\bx, y)  \\
  &= \frac{1}{M^{M/2+1}} \sum_{\bx \in \bXM} \sum_{y \in \cM} \sum_{\substack{\cI \subset [M] \\ |\cI| = k}} \prod_{j \in \cI} \delta(x_j,y) \prod_{\ell \in [M] \setminus \cI} \big(1-\delta(x_\ell,y)\big).
    \end{split}
\end{equation*}
Note that the product on the right-hand side above is
\begin{align*}
  \prod_{j \in \cI} \delta(x_j,y) \prod_{\ell \in [M] \setminus \cI} \big(1-\delta(x_\ell,y)\big)= \delta_{\cI}(\bx,y).
\end{align*}
where $\delta_{\cI}(\bx,y)$ is defined by~\eqref{eqPixyI},
% 
% We let $\Pi(\bx,y; \cI)$ denote the product on the right-hand side above, so that
% \begin{align*}
%   \Pi(\bx,y; \cI):= \prod_{j \in \cI} \delta(x_j,y) \prod_{\ell \in [M] \setminus \cI} \big(1-\delta(x_\ell,y)\big).
% \end{align*}
% Note that, explicitly:
% \begin{equation}\label{eqPixyI}
%   \Pi(\bx,y; \cI) = \begin{cases}
%    1, & \text{if $x_j=y$ for all $j\in\cI$ and 
%    $x_\ell\neq y$ for all $\ell\not\in\cI$;}\\[6pt]
%    0, & \text{else.}
%   \end{cases} 
% \end{equation}

Then, changing the order of summation, we have:
\begin{equation} \label{eq:A(k,M)}
  \begin{split}
  \A(k,M) 
  &= \frac{1}{M^{M/2+1}} \sum_{\substack{\cI \subset [M] \\ |\cI| = k}} \sum_{y \in \cM} \sum_{\bx \in \bXM}   
  \delta_{\cI}(\bx,y)\,.
    \end{split}
\end{equation}

We split the analysis of the sums in~\eqref{eq:A(k,M)} into two cases, depending on whether $\PI$ is empty or not. Accordingly, the expression of the average~\eqref{eq:A(k,M)}
can be conveniently written as 
\begin{align}\label{eqA12}
  \A(k,M) = A_1(k,M) + A_2(k,M), 
\end{align}
where
\begin{align*}%\label{eq:A(k,M)-2}
A_1(k,M)  &= \frac{1}{M^{M/2+1}} 
   \sum_{\substack{\cI \subset [M] \\ |\cI| = k\\ \PI = \emptyset }} 
   \sum_{y \in \cM} 
   \sum_{\bx \in \bXM}   
 \delta_{\cI}(\bx,y) \\[-10pt]
\intertext{and}
A_2(k,M)  &= \frac{1}{M^{M/2+1}} 
  \sum_{\substack{\cI \subset [M] \\ |\cI| = k\\ \PI \neq \emptyset }} 
  \sum_{y \in \cM} 
  \sum_{\bx \in \bXM}   
 \delta_{\cI}(\bx,y). 
\end{align*}

Although there are many subsets $\cI\subset\cM$ in the above sums, the 
involution $\psi$ is the one that sorts them out and allows us to count those that  actually have a non-zero contribution to $A_1(k,M)$ and $A_2(k,M)$.
%%%%%%%%%%%%%%%%%%%%%%%%%%%%
\begin{remark}\label{Remark1}
If $y\in\cI$ is a fixed point of $\psi$ and 
given $\bx\in\bXM$ with a component
$x_j=y$ for some $j\in\cI$, then
  \begin{align*}
     \psi(x_{M+1-j}) = x_j = y = \psi(y) 
  \end{align*}
implies that $x_{M+1-j} = y$.
In view of~\eqref{eqPixyI}, this means that if $\delta_{\cI}(\bx,y)=1$ and $\cP(\cI) = \emptyset$,
then $y$ is not a fixed point of $\psi$.
\end{remark}
%%%%%%%%%%%%%%%%%%%%%%%%%%%%
\begin{remark}\label{Remark2}
Given $\bx\in\bXM$, if $\cP(\cI) \neq \emptyset$
and $y \in \Fpsi$ is a fixed point of $\psi$, then 
$x_j=y$ for some $j\in\cI$, and
  \begin{align*}
     y = \psi(y) = \psi(x_{j}) = x_{M+1-j}.
  \end{align*}
It follows that if $\delta_{\cI}(\bx,y)\neq 0$, then
$x_{M+1-j} \in\cI$, that is, all elements of $\cI$ are in pairs $\{j,{M+1-j}\}$ with both components $x_j$ and $x_{M+1-j}$ in $\Fpsi$.
\end{remark}
Next, we employ these observations to estimate 
$A_1(k,M)$ and $A_2(k,M)$.
% \medskip

%%%%%%%%%%%%%%%%%%%%%%%%%%%%%%%%%%%%%%%%%%%%%%%%%%%%%%%%
\textbf{Case 1.} If $\PI = \emptyset$, by Remark~\ref{Remark1} we know that
the product $\delta_{\cI}(\bx,y)$ is equal to zero if $y \in \Fpsi$. 
For each $y \notin \Fpsi$, a vector $\bx\in\bXM$ 
for which the product 
$\delta_{\cI}(\bx,y)$
is non-zero has $2k$ fixed components.
These are those $k$ components where $x_j=y$ for $j\in\cI$ and,
by the definition of $\bx$ as an element of $\bXM$, the other $k$
components fixed by $\psi$, namely $x_{M+1-j}=\psi(x_j)=\psi(y)$
for $j\in\cI$. 
Therefore, there are 
\begin{align*}
  (M-2)^{M/2 - |\cI|} = (M-2)^{M/2 - k}  
\end{align*}
vectors $\bx$ that have non-zero contribution to $\A(k,M)$. Consequently, we have:
\begin{equation}\label{eqA1cf}
 \begin{split}
    A_1(k,M) &=  \frac{1}{M^{M/2+1}} \sum_{\substack{\cI \subset [M] \\ |\cI| = k\\ \PI = \emptyset }} \sum_{y \in \cM \setminus \Fpsi} (M-2)^{M/2 - k} \\
    &= \frac{|\cM \setminus \Fpsi|}{M^{M/2 + 1}} \binom{M/2}{k} 2^{k} (M-2)^{M/2 - k}.
  \end{split}
\end{equation}
Grouping together factors of similar size and using the assumption on the estimate of the number of fixed points of $\psi$, we have:
\begin{equation}\label{eqA1estimate}
  \begin{split}
    A_1(k,M) 
    &= \frac{|\cM \setminus \Fpsi|}{ M} \cdot \dfrac{(M/2) (M/2 - 1) \cdots (M/2 - k + 1)}{k! (M/2)^k} \left(1 - \frac{2}{M}\right)^{M/2 - k}\\
    &= \frac{1-\theta }{k!}\left(1+O\big(M^{\eta-1}\big)\right)
    \cdot \left(1+O\big(k^2M^{-1}\big)\right)
    \cdot e^{-1}\left(1+O\big(kM^{-1}\big)\right)\\
     &= \frac{1-\theta}{e\,k!}
     \left(1+O\big(M^{\eta-1}+ k^2M^{-1}\big)\right)
  \end{split}
\end{equation}
for $M$ tending to infinity. 
% \medskip

%%%%%%%%%%%%%%%%%%%%%%%%%%%%%%%%%%%%%%%%%%%%%%%%%%%%%%%%
\textbf{Case 2.} 
If $\PI \neq \emptyset$, 
following Remark~\ref{Remark2}, we know that \mbox{$\delta_{\cI}(\bx,y)=0$} 
unless $y \in \Fpsi$ and 
$\{j, M+1 - j\} \subset \cI$ for all $j \in \cI$. 
Since $M$ is even,
this case only happens when $|\cI| = k$ is also even, which means $A_2(k,M) = 0$ when $k$ is odd. When $k$ is even, for such $y$ and $\cI$, there are 
\begin{align*}
  (M-1)^{M/2 - |\cI|/2} = (M-1)^{M/2 - k/2}
  \end{align*}
    vectors $\bx$ for which $\delta_\cI(\bx,y) = 1$. 
Then, since the subsets $\cI$ that have a non-zero contribution to
$A_2(k,M)$ have the property that
$\{j, M+1 - j\} \subset \cI$ for all $j \in \cI$,
it means that they also satisfy \mbox{$|\PI| = k/2$}; therefore the sum $A_2(k,M)$ reduces to
\begin{equation}\label{eqA2cf}
  \begin{split}
 A_2(k,M) &=   \frac{1}{M^{M/2+1}} \sum_{\substack{\cI \subset [M] \\ |\cI| = k\\ |\PI| = k/2 }} \sum_{y \in  \Fpsi} (M-1)^{M/2 - k/2} \\
 &= \frac{|\Fpsi|}{M^{M/2+1}} \binom{M/2}{k/2} (M-1)^{M/2 - k/2} .
\end{split}
\end{equation}
Grouping the factors conveniently, we obtain:
\begin{equation}\label{eqA2estimate}
  \begin{split}
  A_2(k,M) &= \frac{|\Fpsi|}{M} 
  \cdot \frac{(M/2)(M/2-1)\cdots (M/2 - k/2+1)}{(k/2)!M^{k/2}}
  \left(1 - \frac{1}{M}\right)^{M/2 - k/2}\\
 &= \frac{\theta}{(k/2)!} \left(1 +O\big(M^{\eta-1}\big)\right)
    \cdot  2^{-k/2} \left(1+O\big(k^2M^{-1}\big)\right)
    e^{-1/2}\left(1+O\big(kM^{-1}\big)\right)\\
   &= \frac{\theta}{e^{1/2}\, 2^{k/2}(k/2)!}  
    \left(1+O\big(M^{\eta-1}+ k^2M^{-1}\big)\right)
  \end{split}   
\end{equation}
as $M$ tends to infinity.

Bringing together relations~\eqref{eqA1cf} and~\eqref{eqA2cf} into~\eqref{eqA12} 
and taking into account that $A_2(k,M)= 0$ for odd $k$, 
we derive the closed-form expression of~$A(k,M)$.
\begin{lemma}\label{LemmaAverageCF}
The average 
$\A(k, M) 
    = \frac{1}{|\bXM|}  \sum_{ \bx \in \bXM} {m_k(\bx)}{/M}$
equals:
    % We have:
\begin{equation*}%\label{eqLemmaAverageCF}
  \begin{split}
%    &\ \A(k, M) 
%     = \frac{1}{|\bXM|}  \sum_{ \bx \in \bXM} \frac{m_k(\bx)}{M}\\[6pt]
%     &=
&
     \begin{cases}
   \begin{aligned}
     \frac{|\cM \setminus \Fpsi|}{M^{M/2 + 1}} \binom{M/2}{k} 2^{k} (M-2)^{M/2 - k}, 
    \end{aligned} 
     & \text{if $k$ is odd}; \\[18pt]
   \begin{aligned}
     \frac{|\cM \setminus \Fpsi|}{M^{M/2 + 1}} 
     &\binom{M/2}{k} 2^{k} (M-2)^{M/2 - k}\\[6pt]
     & + \frac{|\Fpsi|}{M^{M/2+1}} \binom{M/2}{k/2} (M-1)^{M/2 - k/2}.
   \end{aligned}
     & \text{if $k$ is even;}
   \end{cases} 
   \\[10pt]
   %%%%%%%%%%%%%%%%%%%%%%%%%%%%%%%
%    &=
 =&
     \begin{cases}
   \begin{aligned}
     \frac{|\cM \setminus \Fpsi|}{M} \cdot \frac{1}{ek!}
     \big(1+O(k^2M^{-1}) \big), 
    \end{aligned} 
     & \text{if $k$ is odd}; \\[16pt]
   \begin{aligned}
     \bigg(
      \frac{|\cM \setminus \Fpsi|}{M} \cdot \frac{1}{ek!}
      +
      \frac{|\Fpsi|}{M} \cdot \frac{1}{e^{1/2}2^{k/2}(k/2)!}
     \bigg)
      \big(1+O(k^2M^{-1}) \big), 
   \end{aligned}
     & \text{if $k$ is even.}
   \end{cases}    
  \end{split}
\end{equation*}
\end{lemma}
Next, under our assumption~\eqref{eqNumberOfFixedPoints} 
regarding the number of fixed points of~$\psi$, 
by inserting the estimates~\eqref{eqA1estimate} and~\eqref{eqA2estimate} into~\eqref{eqA12},
we derive the following estimate of the average.
% We state in the next lemma the result obtained 
% by inserting the estimates~\eqref{eqA1estimate} and~\eqref{eqA2estimate} into~\eqref{eqA12}.

%%%%%%%%%%%%%%%%%%%%%
\begin{lemma}\label{LemmaAverage}
Let $\theta\in [0,1]$ and $\eta \in [0,1)$ and suppose
that the number of fixed points of~$\psi$~is 
$\theta M\left(1 + O(M^{\eta-1})\right)$ as $M$ tends to infinity.
Then: 
% we~have:
\begin{equation*}%\label{eqLemmaAverage}
  \begin{split}
   \A(k, M) 
%     &= \frac{1}{|\bXM|}  \sum_{ \bx \in \bXM} \frac{m_k(\bx)}{M}\\
    &=
     \begin{cases}
     \fdfrac{1-\theta}{ek!}
     \left(1+O\big(M^{\eta-1}+ k^2M^{-1}\big)\right), & \text{if $k$ is odd}; \\[6pt]
     \left(\fdfrac{1-\theta}{e\,k!}
     + \fdfrac{\theta}{e^{1/2}\, 2^{k/2}(k/2)!}  \right)
    \left(1+O\big(M^{\eta-1}+ k^2 M^{-1}\big)\right), & \text{if $k$ is even.}
   \end{cases}   
  \end{split}
\end{equation*}

\end{lemma}

%%%%%%%%%%%%%%%%%%%%%
\section{The second moment about the mean}\label{SectionMoment}
Since the average depends on the parity of $k$, we need to separate the discussion into two cases.

%%%%%%%%%%%%%%%%%%%%%%%%%%%%%%%%%%%%%%%%%%%%%%%%%%%%%%%%%
\subsection{\texorpdfstring{The case $k$ is odd}{The case k is odd}}

Considering the main term in the formula from Lemma~\ref{LemmaAverage} when~$k$ is odd, we let:
\begin{align}\label{eqM2odd0}
    \M(k,M) := \frac{1}{|\bXM|} \sum_{\bx \in \bXM} \left(\frac{m_k(\bx)}{M} - \frac{1-\theta}{ek!}\right)^2,\ \ \text{if $k$ is odd.} 
\end{align} 

Expanding the binomial, we find that $\M(k,M)$ can be written as 
\begin{align}\label{eqM2odd}
    \M(k,M) = \left(\frac{1-\theta}{ek!}\right)^2 -  \frac{2(1-\theta)}{ek!} \A(k, M) + \frac{1}{M^{M/2+2}} S_2(k,M),
\end{align}
where
\begin{align}\label{eqS2def}
    S_2(k,M) := %frac{1}{M^{M/2+2}} 
    \sum_{\bx \in \bXM} m_k(\bx)^2.
\end{align}
Our aim is to show that $\M(k,M)$ is small as $M \to \infty$
and, for this, we will show that
the main term of $S_2(k,M)$ equals the square of the main term 
of $\A(k,M)$ given in Lemma~\ref{LemmaAverage}.

By \eqref{eq:m_k(x)}, we see that 
\begin{align*}
    S_2(k, M) &= %\frac{1}{M^{M/2+2}} 
    \sum_{\bx \in \bXM} 
    \left(\sum_{y \in \cM} \delta_k(\bx,y)\right)^2 
    = %\frac{1}{M^{M/2 + 2}} 
    \sum_{\bx \in \bXM} 
    \sum_{y \in \cM}\sum_{y' \in \cM} \delta_k(\bx, y) \delta_k(\bx, y').
\end{align*}
Separating the sums of the diagonal and the off-diagonal yields:
\begin{equation}\label{eqS2Def}
 \begin{split}
    S_2(k, M)
    &= %\frac{1}{M^{M/2 + 2}} 
    \sum_{\bx \in \bXM} 
    \sum_{y \in \cM} \sum_{\substack{y' \in \cM\\y'= y}}
    \delta_k(\bx, y) \delta_k(\bx, y')
    + %\frac{1}{M^{M/2 + 2}} 
    \sum_{\bx \in \bXM} 
    \sum_{y \in \cM} \sum_{\substack{y' \in \cM\\y'\neq y}}
     \delta_k(\bx, y) \delta_k(\bx, y') \\
    &=  
     M^{M/2+1} \A(k, M)
     + U(k,M),
\end{split}
\end{equation}
where
\begin{align*}
    U(k,M):&=  \sum_{\bx \in \bXM} 
    \sum_{y \in \cM} \sum_{\substack{y' \in \cM\\y'\neq y}}
     \delta_k(\bx, y) \delta_k(\bx, y').
\end{align*}
This can be written as:
\begin{align*}
    U(k,M) 
    &= 
    \sum_{\bx \in \bXM} \sum_{y \in \cM} \sum_{\substack{y' \in \cM\\y'\neq y}} \sum_{\substack{\cI, \cI' \subset [M] \\ |\cI| = |\cI'| = k}} \delta_{\cI}(\bx, y) \delta_{\cI'}(\bx, y') \notag \\
    &= 
    \sum_{\substack{\cI, \cI' \subset [M] \\ |\cI| = |\cI'| = k}} \sum_{y \in \cM} \sum_{\substack{y' \in \cM\\y'\neq y}}  \sum_{\bx \in \bXM}  \delta_{\cI}(\bx, y) \delta_{\cI'}(\bx, y') .
\end{align*}

By the argument in Remark~\ref{Remark2}, for odd $k$, 
$\delta_\cI(\bx, y)$ is zero unless $\cI$ contains no pair $\{j, M+1 - j\}$ for every $j \in [M]$. 
Also, since $y \neq y'$, the product
$\delta_{\cI}(\bx, y) \delta_{\cI'}(\bx, y')$ is equal to zero if $\cI \cap \cI' \neq \emptyset$. Therefore, we only need to consider the case $\PI = \PIp = \emptyset$ and $\cI \cap \cI' = \emptyset$ for estimating $U(k,M)$. 
In other words,
\begin{align}\label{eqDefU}
    U(k,M) = %\frac{1}{M^{M/2 + 2}}
    \sum_{\substack{\cI, \cI' \subset [M] \\ |\cI| = |\cI'| = k \\ \cI \cap\cI'= \emptyset \\
    \PI = \PIp = \emptyset}} \sum_{y \in \cM} \sum_{\substack{y' \in \cM\\y'\neq y}}  \sum_{\bx \in \bXM}  \delta_{\cI}(\bx, y) \delta_{\cI'}(\bx, y').
\end{align}
We further separate the terms in the sum depending on the 
size of the set $\cP(\cI\cup\cI')$:
\begin{align}\label{eqDefUkMU1U2}
    U(k,M) &= U_1(k,M) + U_2(k,M),
\end{align}
where
\begin{align*}
    U_1(k,M) &= %\frac{1}{M^{M/2 + 2}}
    \sum_{\substack{\cI, \cI' \subset [M] \\ |\cI| = |\cI'| = k \\ \cI \cap\cI'= \emptyset \\
    \PI = \PIp = \emptyset\\
    \PIIp \neq \emptyset
    }} \sum_{y \in \cM} \sum_{\substack{y' \in \cM\\y'\neq y}}  \sum_{\bx \in \bXM}  \delta_{\cI}(\bx, y) \delta_{\cI'}(\bx, y')\\
\intertext{and}
    U_2(k,M) &= %\frac{1}{M^{M/2 + 2}}
    \sum_{\substack{\cI, \cI' \subset [M] \\ |\cI| = |\cI'| = k \\ \cI \cap\cI'= \emptyset \\
    \PI = \PIp = \emptyset\\
    \PIIp = \emptyset
    }} \sum_{y \in \cM} \sum_{\substack{y' \in \cM\\y'\neq y}}  \sum_{\bx \in \bXM}  \delta_{\cI}(\bx, y) \delta_{\cI'}(\bx, y').
\end{align*}

\begin{remark}\label{Remark3}
Following conditions under the summations in $U_1(k,M)$, let us suppose that $\cI,\cI'\subset [M]$ and $\PIIp \neq \emptyset$.
\begin{enumerate}
    \item[(a)]
 Then, if $y\neq y'$ and $\bx\in \bXM$ are such that
$\delta_{\cI}(\bx, y) \delta_{\cI'}(\bx, y')=1$, 
it follows that there exists 
$j\in\cI$ such that 
$M+1-j\in\cI'$,    
$x_j=y$, and 
$x_{M+1-j}=y'$.
Then 
\begin{align*}
    \psi(y) = \psi(x_j) = x_{M+1-j} = y',
\end{align*}
which implies $\psi(y')=y$ and $\psi(y)=y'$.
Additionally, it should be noted that it follows that $\cI' = \{M+1-j: j\in\cI \}$.
    \item[(b)]
Moreover, under the same conditions, if there is another component, say $x_h$ with $h\in\cI$,
in order to still have $\delta_{\cI}(\bx, y)\neq 0$,
then $x_h=y$, so that $\psi(x_h)=y'$.
\end{enumerate}
Then, by combining (a) and (b), we find that 
$\PIIp$ must have the same cardinality as both 
$\cI$ and $\cI'$.
\end{remark}
%%%%%%%%%%%%%%%%%%%%%%%%%%%%%%%%%%%%%%%%%%%%%%%%%%%%%%%%
\textbf{Case 1.} 
To estimate $U_1(k,M)$, the condition that needs to be satisfied is $\PIIp \neq \emptyset$.
Then, based on Remark~\ref{Remark3}, a product 
$\delta_{\cI}(\bx, y) \delta_{\cI'}(\bx, y')$ 
does contribute to $U_1(k,M)$ when 
$|\PIIp|=k$ and the vector $\bx$ has exactly $2k$ components fixed, namely $x_j = y$ for all~$j \in \cI$ and $x_{\ell} = y'$ for all $\ell \in \cI'$. 
It follows that there are $(M-2)^{M/2 - k}$
such vectors $\bx$, that is,
\begin{align*}
  \sum_{\bx \in \bXM} \delta_{\cI}(\bx,y) 
  \delta_{\cI'}(\bx, y') = (M-2)^{M/2 - k}.
\end{align*}
This holds for any $y \neq y' \in \cM$ with $y' = \psi(y)$,
meaning that $y$ is not a fixed point of $\psi$, and neither is $y'$.  Therefore 
\begin{align*}
    U_1(k,M) &= %\frac{1}{M^{M/2 + 2}}
    \sum_{\substack{\cI, \cI' \subset [M] \\ |\cI| = |\cI'| = k \\ \cI \cap \cI'= \emptyset \\
    \cP(\cI)=\cP(\cI')=\emptyset\\
    |\PIIp| = k}} \sum_{y \in \cM} 
    \sum_{\substack{y' \in \cM\\y'\neq y\\y'=\psi(y)}}  
    (M-2)^{M/2 - k} \\
    &= %\frac{|\cM \setminus \Fpsi|}{M^{M/2+2}} 
    |\cM \setminus \Fpsi| %\cdot
    \binom{M/2}{k} 2^k (M-2)^{M/2 - k}. \\ 
    % &= \frac{|\cM \setminus \Fpsi|}{M^2 k!} \cdot \frac{M/2(M/2-1) \cdots (M/2 - k + 1)}{(M/2)^k} \left(1 - \frac{2}{M}\right)^{M/2 - k},
\end{align*}
By appropriately grouping the factors of the product, we find that
\begin{equation}\label{eqU1}
  \begin{split}
   \frac{U_1(k,M)}{M^{M/2 + 2}} &=  
    \frac{|\cM \setminus \Fpsi|}{k! M^2} \cdot \frac{(M/2)(M/2-1) \cdots (M/2 - k + 1)}{(M/2)^k} \left(1 - \frac{2}{M}\right)^{M/2 - k}\\
    &=\frac{1-\theta}{k!M}
    \left(1 +O\big(M^{\eta-1}\big)\right)
     \left(1+O\big(k^2M^{-1}\big)\right)
    e^{-1}\left(1+O\big(kM^{-1}\big)\right)\\
   &= \frac{\theta}{e\, k!}  \cdot \frac{1}{M}
    \left(1+O\big(M^{\eta-1}+ k^2M^{-1}\big)\right),
  \end{split}
\end{equation}
which is $O(M^{-1})$ when $k$ is fixed and $M \to \infty$. 

%%%%%%%%%%%%%%%%%%%%%%%%%%%%%%%%%%%%%%%%%%%%%%%%%%%%%%%%
\textbf{Case 2.} 
If $\PIIp = \emptyset$, then $\delta_{\cI}(\bx, y) \delta_{\cI'}(\bx, y')$ is zero unless 
\mbox{$y, y' \notin \Fpsi$} and $y \neq \psi(y')$. Therefore, each vector $\bx$ for which $\delta_{\cI}(\bx, y) \delta_{\cI'}(\bx, y')$ is non-zero has $4k$ components fixed, namely 
$x_j = y$, $ x_{M+1 - j} = \psi(y)$ for all $j \in \cI$ 
and $x_\ell = y'$, $x_{M+1 - \ell} = \psi(y')$ 
for all $\ell \in \cI'$. 
This shows that 
\begin{align*}
\sum_{\bx \in \bXM} \delta_{\cI}(\bx,y) \delta_{\cI'}(\bx, y') = (M-4)^{M/2 - 2k}.\end{align*}
This holds for any $y \neq y' \in \cM \setminus \Fpsi$ with $y \neq \psi(y')$. It follows, 
\begin{align*}
    U_2(k,M)  
    &= %\frac{1}{M^{M/2 + 2}} 
    \sum_{\substack{y \neq y' \in \cM \setminus \Fpsi \\ y \neq \psi(y')}} \sum_{\substack{\cI, \cI' \subset [M] \\ |\cI| = |\cI'| = k \\ 
    \cI \cap\cI'= \emptyset \\
    \PI = \PIp = \emptyset\\
    \PIIp = \emptyset}}  (M-4)^{M/2 - 2k} \\
    &= 
    |\cM \setminus \Fpsi|\left(|\cM \setminus \Fpsi| - 2\right)
    \binom{M/2}{k} 2^k \binom{M/2-k}{k} 2^k (M-4)^{M/2 - 2k}.
\end{align*}
To estimate $U_2(k,M)$ we group the factors as follows:
\begin{equation}\label{eqU2}
  \begin{split}
   \frac{U_2(k,M)}{M^{M/2 + 2}} &=
   \frac{|\cM \setminus \Fpsi|\left(|\cM \setminus \Fpsi| - 2\right)}{M^{2}  } \cdot  \frac{(M/2) \cdots (M/2 - k + 1)}{k! (M/2)^k}\\
   &\quad
   \times \frac{(M/2 - k) \cdots (M/2 - 2k +1)}{k! (M/2)^k} \left(1 - \frac{4}{M}\right)^{M/2 - 2k}\\
   &=\frac{(1-\theta)^2}{(k!)^2}
    \left(1 +O\big(M^{\eta-1}\big)\right)
    \left(1+O\big(k^2M^{-1}\big)\right)
    e^{-2}\left(1+O\big(kM^{-1}\big)\right)\\
     &=\left(\frac{1-\theta}{e\,k!}\right)^2
      \left(1+O\big(M^{\eta-1}+ k^2M^{-1}\big)\right).
  \end{split}
\end{equation}

On combining the estimates~\eqref{eqU1},~\eqref{eqU2}, and~\eqref{eqDefUkMU1U2}, it yields
\begin{equation}\label{eqEstimateU}
  \begin{split}
   \frac{U(k,M)}{M^{M/2 + 2}} &= 
   \left(\frac{1-\theta}{e\,k!}\right)^2
      \left(1+O\big(M^{\eta-1}+ k^2M^{-1}\big)\right).
  \end{split}
\end{equation}
% as $M$ tends to infinity.
Next, knowing from Lemma~\ref{LemmaAverage} that 
$\A(k,M) = O(1)$ and inserting the estimate~\eqref{eqEstimateU}
into~\eqref{eqS2Def},
we also obtain:
\begin{equation}\label{eqS2kodd}
  \begin{split}
  \frac{S_2(k, M)}{M^{M/2 + 2}} 
  &=\frac{\A(k,M)}{M} + \frac{U(k,M)}{M^{M/2 + 2}}\\
  &=\left(\frac{1-\theta}{e\,k!}\right)^2
      \left(1+O\big(M^{\eta-1}+ k^2M^{-1}\big)\right).
  \end{split}
\end{equation}
Then, on using~\eqref{eqS2kodd} and the estimate from Lemma~\ref{LemmaAverage} 
for $k$ odd in~\eqref{eqM2odd}, it follows
\begin{equation}\label{eqM2koddfinal}
  \begin{split}
  \M(k,M) &= 
   \left(\frac{1-\theta}{ek!}\right)^2 
   -  \frac{2(1-\theta)}{ek!} \left(\frac{1-\theta}{ek!}\right)
   \left(1+O\big(M^{\eta-1}+ k^2M^{-1}\big)\right)\\
    &\quad
   + \left(\frac{1-\theta}{e\,k!}\right)^2
      \left(1+O\big(M^{\eta-1}+ k^2M^{-1}\big)\right)\\
    &=  O\big(M^{\eta-1}+ k^2M^{-1}\big)
  \end{split}
\end{equation}
as $k$ is fixed and $M \to \infty$, for $k$ odd.

%%%%%%%%%%%%%%%%%%%%%%%%%%%
\subsection{\texorpdfstring{The case $k$  even}{The case k even}}
Using the main term in the formula from Lemma~\ref{LemmaAverage} when $k$ is even, in this case, for the second moment about the mean, we let:
\begin{align}\label{eqM2even}
    \M(k,M) := \frac{1}{|\bXM|} \sum_{\bx \in \bXM} \left(\frac{m_k(\bx)}{M} - \frac{1-\theta}{ek!}
    - \frac{\theta}{e^{1/2}\,2^{k/2} (k/2)! }
    \right)^2,\ \ \text{if $k$ is even.} 
\end{align} 
Expanding the square and denoting, as in the previous case, by $S_2(k,M)$
the sum of the $k$-counter functions, a notation that does not depend on the parity of $k$, we have
\begin{align*}
  \M(k,M) &= 
  \left(\frac{1-\theta}{e\,k!} + \frac{\theta}{e^{1/2}\, 2^{k/2} (k/2)!}\right)^2 \\
  &\quad- 2\left(\frac{1-\theta}{e\,k!} + \frac{\theta}{e^{1/2}\, 2^{k/2} (k/2)!}\right) \A(k,M) 
  + \frac{S_2(k,m)}{M^{M/2+2}}.
\end{align*}
After similar computation as in the case $k$ is odd, we obtain 
\begin{equation}\label{eqS2DefV}
 \begin{split}
    S_2(k, M) =
     M^{M/2+1} \A(k, M)
     + V(k,M),
    %\frac{1}{M^{M/2 + 2}} 
\end{split}
\end{equation}
where $V(k,M)$ verifies the same formula as the previous $U(k,M)$:
\begin{align*}
    V(k,M):&=  
    \sum_{\substack{\cI, \cI' \subset [M] \\ |\cI| = |\cI'| = k}} \sum_{y \in \cM} \sum_{\substack{y' \in \cM\\y'\neq y}}  \sum_{\bx \in \bXM}  \delta_{\cI}(\bx, y) \delta_{\cI'}(\bx, y') .
\end{align*}
Following that, unlike the case when $k$ was odd,
we split $V(k,M)$ into shorter pieces as follows: 
$V_1(k,M)$ is the sub-sum when \mbox{$\PI = \PIp = \emptyset$}, 
$V_2(k,M)$ is the sub-sum when either \mbox{$\PI \neq \emptyset$} 
or $\PIp \neq \emptyset$ but not both, and 
$V_3(k,M)$ is the sub-sum when both $\PI$ and~$\PIp$ 
are non-empty. 
Again, since $y \neq y'$,  the product 
$\delta_{\cI}(\bx, y) \delta_{\cI'}(\bx, y')$ is equal to zero if $\cI \cap \cI' \neq \emptyset$. Thus
\begin{align}\label{eqV}
    V(k,M)  = V_1(k,M) + V_2(k, M) + V_3(k,M),
\end{align}
where
\begin{align*}
    V_1(k,M) &=  %\frac{1}{M^{M/2 + 2}}
    \sum_{\substack{\cI, \cI' \subset [M] \\ |\cI| = |\cI'| = k \\
    \cI \cap \cI' = \emptyset \\\PI = \PIp = \emptyset}} \sum_{y \in \cM} \sum_{\substack{y' \in \cM\\y'\neq y}}  \sum_{\bx \in \bXM}  \delta_{\cI}(\bx, y) \delta_{\cI'}(\bx, y'), \\
V_2(k,M) &=  %\frac{1}{M^{M/2 + 2}}
   2\times \sum_{\substack{\cI, \cI' \subset [M] \\ |\cI| = |\cI'| = k \\ \cI \cap \cI' = \emptyset \\\PI \neq \emptyset,\; \PIp = \emptyset}} \sum_{y \in \cM} \sum_{\substack{y' \in \cM\\y'\neq y}}  \sum_{\bx \in \bXM}  \delta_{\cI}(\bx, y) \delta_{\cI'}(\bx, y'), \\
\intertext{(Here, the coefficient $2$ accounts for the other sum where $\PI = \emptyset,$\; $\PIp \neq \emptyset$, a sum that, because of the symmetry, is equal to the one shown.)}
V_3(k,M) &= % \frac{1}{M^{M/2 + 2}}
\sum_{\substack{\cI, \cI' \subset [M] \\ |\cI| = |\cI'| = k \\ \cI \cap \cI' = \emptyset \\\PI \neq \emptyset,\; \PIp \neq \emptyset}} \sum_{y \in \cM} \sum_{\substack{y' \in \cM\\y'\neq y}}  \sum_{\bx \in \bXM}  \delta_{\cI}(\bx, y) \delta_{\cI'}(\bx, y').
\end{align*}

Next, we proceed to examine each sum $V_1(k,M)$, $V_2(k,M)$ and $V_3(k, M)$ separately.
%%%%%%%%%%%%%%%%%%%%%%%%%%%%%%%%%%%%%%%%%%%%%%%%%%%%%%%
\subsubsection{\texorpdfstring{Estimating $V_1(k,M)$}{Estimating V1(k,M)}}
The sum $V_1(k,M)$ is identical to $U(k,M)$ from~\eqref{eqDefU}, and its estimate~\eqref{eqEstimateU} does not depend on the parity of $k$. Hence 
\begin{equation}\label{eqV1}
  \begin{split}
  % V_1(k,M) 
  \frac{V_1(k,M)}{M^{M/2 + 2}}
  =\left(\frac{1-\theta}{e\,k!}\right)^2
      \left(1+O\big(M^{\eta-1}+ k^2M^{-1}\big)\right)
  \end{split}
\end{equation}
for $M \to \infty$.
%%%%%%%%%%%%%%%%%%%%%%%%%%%%%%%%%%%%%%%%%%%%%%%%%%%%%%%
\subsubsection{\texorpdfstring{Estimating $V_2(k,M)$}{Estimating V2(k,M)}}
As $\PI \neq \emptyset$, $\delta_\cI(\bx,y)$ is zero unless 
\mbox{$\{j, M+1 - j\} \subset \cI$} for all $j \in \cI$ and $y \in \Fpsi$. On the other hand, as $\PIp = \emptyset$, the necessary condition for which $\delta_{\cI'}(\bx, y')  = 1$ is $y' \notin \Fpsi$. 
Now, for any $y \in \Fpsi, y' \notin \Fpsi$, there are 
\begin{align*}
  (M-3)^{M/2 - |\cI|/2 - |\cI'|} = (M-3)^{M/2 - 3k/2}
\end{align*}
vectors $\bx$ that satisfy $\delta_{\cI}(\bx, y) \delta_{\cI'}(\bx, y')   = 1$, since each such vector $\bx$ has~$3k$ components fixed while the other components must be different from 
\mbox{$y=\psi(y)$}, $y'$, and $\psi(y')$, together with the relation $x_j = \psi(x_{M+1 - j})$ for all $j \in [M]$. Thus 
\begin{align*}
  V_2(k,M) &=  
  2|\Fpsi| |\cM \setminus \Fpsi|\times
  \sum_{\substack{\cI, \cI' \subset [M] \\ |\cI| = |\cI'| = k \\ \cI \cap \cI' = \emptyset \\\PI \neq \emptyset,\; \PIp = \emptyset}} (M-3)^{M/2 - 3k/2}.
\end{align*}
This further equals
\begin{align*}
  V_2(k,M)
  &=  2|\Fpsi| |\cM \setminus \Fpsi| \binom{M/2}{k/2} \binom{M/2 - k/2}{k} 2^{k} (M-3)^{M/2 - 3k/2}.
\end{align*}
Next, we group factors appropriately to find that
\begin{equation}\label{eqV2}
  \begin{split}
  \frac{V_2(k,M)}{M^{M/2 + 2}}
  % &=  \frac{2|\Fpsi| |\cM \setminus \Fpsi|}{M^2} \binom{M/2}{k/2} \binom{M/2 - k/2}{k} 2^{k} (M-3)^{M/2 - 3k/2}\\
  &=  \frac{2|\Fpsi| |\cM \setminus \Fpsi|}{M^2} \cdot \frac{\prod_{j=0}^{k/2-1} (M/2 - j)}{(k/2)!M^k} \\
  &\phantom{
  =  \frac{2|\Fpsi| |\cM \setminus \Fpsi|}{M^2} \cdot 
  }
  \times \frac{\prod_{\ell=0}^{k-1}(M/2-k/2-\ell)}{k!(M/2)^k} \left(1 - \frac{3}{M}\right)^{M/2 - 3k/2}\\
  &= \frac{2\theta(1-\theta)}{(k/2)!k!}
    \left(1 +O\big(M^{\eta-1}\big)\right)
    \cdot  \left(1+O\big(k^2M^{-1}\big)\right)
    e^{-3/2}\left(1+O\big(kM^{-1}\big)\right)\\
   &= \frac{2\theta(1-\theta)}{e^{3/2}(k/2)!k!}
    \left(1+O\big(M^{\eta-1}+ k^2M^{-1}\big)\right).
  \end{split}
\end{equation}

%%%%%%%%%%%%%%%%%%%%%%%%%%%%%%%%%%%%%%%%%%%%%%%%%%%%%%%%
\subsubsection{\texorpdfstring{Estimating $V_3(k,M)$}{Estimating V3(k,M)}}
Similarly with above arguments, the product 
$\delta_{\cI}(\bx, y) \delta_{\cI'}(\bx, y')$ 
is zero unless 
\mbox{$\{j, M+1 -j\} \subset \cI$} for all $j \in \cI$, 
\mbox{$\{\ell, M+1 -\ell\} \subset \cI'$} for all $\ell \in \cI'$, and $y, y' \in \Fpsi$. For such $\cI, \cI'$, there are
\begin{align*}
(M-2)^{M/2 - |\cI|/2 - |\cI'|/2} = (M-2)^{M/2 - k}\end{align*}
vectors $\bx$ that satisfy $\delta_{\cI}(\bx, y) \delta_{\cI'}(\bx, y') = 1$. This holds for all $y, y' \in \Fpsi$ with $y \neq y'$. It follows, 
\begin{align*}
    V_3(k, M) &= 
    % \frac{|\Fpsi|(|\Fpsi| -1)}{M^{M/2 + 2}} 
    |\Fpsi|(|\Fpsi| -1)\times
    \sum_{\substack{\cI, \cI' \subset [M] \\ |\cI| = |\cI'| = k \\ \cI \cap \cI' = \emptyset \\\PI \neq \emptyset,\; \PIp \neq \emptyset}} (M-2)^{M/2 - k}\\
     &= 
    |\Fpsi|(|\Fpsi| -1)\times
    \binom{M/2}{k/2}  \binom{M/2-k/2}{k/2}
    (M-2)^{M/2 - k}.
\end{align*}
Grouping the factors in order to find the asymptotic behavior, we obtain
\begin{equation}\label{eqV3}
  \begin{split}
    \frac{V_3(k,M)}{M^{M/2 + 2}}
    &= \frac{|\Fpsi|(|\Fpsi| -1)}{M^{2}} \binom{M/2}{k/2} M^{-k/2} \binom{M/2-k/2}{k/2} M^{-k/2} \left(1 - \frac{2}{M}\right)^{M/2 - k}\\
  &= \frac{\theta^2}{(k/2)!(k/2)!}
    \left(1 +O\big(M^{\eta-1}\big)\right)
     \left(1+O\big(k^2M^{-1}\big)\right)
    e^{-1}\left(1+O\big(kM^{-1}\big)\right)\\
   &= \left(\frac{\theta}{e^{1/2}(k/2)!}\right)^2
    \left(1+O\big(M^{\eta-1}+ k^2M^{-1}\big)\right).
  \end{split}
\end{equation}

%%%%%%%%%%%%%%%%%%%%%%%%%%%%%%%%%%%%%%%%%%%%%%%%%%%%%%%%
\subsubsection{\texorpdfstring{Ending the estimation of the  second moment when $k$ is even}{Ending the estimation of the  second moment when k is even}}

On inserting the estimates~\eqref{eqV3}, \eqref{eqV2}, 
\eqref{eqV1} into~\eqref{eqV}, then using the result 
into~\eqref{eqS2DefV}, and knowing from Lemma~\ref{LemmaAverage} that 
% \mbox{$\A(k,M) = O(1)$}, 
\mbox{$\A(k,M)/M = O(M^{-1})$}, 
it yields:
\begin{equation}\label{eqS2keven}
  \begin{split}
  \frac{S_2(k, M)}{M^{M/2 + 2}} 
  &=\frac{\A(k,M)}{M}+\frac{V(k, M)}{M^{M/2 + 2}} \\
  &= %O(M^{-1})+ 
  \left(\left(\sdfrac{1-\theta}{e\,k!}\right)^2 + \sdfrac{2\theta(1-\theta)}{e^{3/2}(k/2)! k!} + \left(\sdfrac{\theta }{e^{1/2}\,(k/2)! }\right)^2\right)
      \left(1+O\big(M^{\eta-1}+ k^2M^{-1}\big)\right)\\
  &= \left(\frac{1-\theta}{e\,k!} + \frac{\theta}{e^{1/2}\, 2^{k/2} (k/2)!}\right)^2
   \left(1+O\big(M^{\eta-1}+ k^2M^{-1}\big)\right).
  \end{split}
\end{equation}

Then, inserting~\eqref{eqS2keven} and the estimate from 
Lemma~\ref{LemmaAverage} 
for $k$ even into~\eqref{eqM2even}, it follows
\begin{equation}\label{eqM2kevenfinal}
  \begin{split}
  \M(k,M) &= 
  \left(\frac{1-\theta}{e\,k!} + \frac{\theta}{e^{1/2}\, 2^{k/2} (k/2)!}\right)^2 \\
  &\quad - 2\left(\frac{1-\theta}{e\,k!} + \frac{\theta}{e^{1/2}\, 2^{k/2} (k/2)!}\right) ^2
    \left(1+O\big(M^{\eta-1}+ k^2 M^{-1}\big)\right)
  \\
  &\quad
  % + \frac{S_2(k,m)}{M^{M/2+2}}.
  +\left(\frac{1-\theta}{e\,k!} + \frac{\theta}{e^{1/2}\, 2^{k/2} (k/2)!}\right)^2
   \left(1+O\big(M^{\eta-1}+ k^2M^{-1}\big)\right)\\
   &= O\big(M^{\eta-1}+ k^2M^{-1}\big).
  \end{split}
\end{equation}
% as $k$ is fixed and $M \to \infty$, for $k$ even.

To sum up, we state in the following lemma 
the estimate that holds,
according to~\eqref{eqM2koddfinal} and~\eqref{eqM2kevenfinal},
in both cases of $k$ being odd or even.

%%%%%%%%%%%%%%%%%%%%%%%%%%%%%%%%%%%%%%
\begin{lemma}\label{LemmaMoment}
We have:
\begin{equation}\label{eqLemmaMoment}
  \begin{split}
   \M(k, M) 
    &= O\big(M^{\eta-1}+ k^2 M^{-1}\big).
  \end{split}
\end{equation}
\end{lemma}

%%%%%%%%%%%%%%%%%%%%%%%%%%%%%%%%%%%%%%%%
\section{Proof of Theorem~\ref{Theorem1}}\label{SectionProofTheorem1}
We let $MT(k)$ be the shorthand notation for the main term of the estimate of the average from Lemma~\ref{LemmaAverage}:
\begin{align*}
    MT(k) = \begin{cases}
     \fdfrac{1-\theta}{ek!}, & \text{if $k$ is odd}; \\[6pt]
     \fdfrac{1-\theta}{e\,k!}
     + \fdfrac{\theta}{e^{1/2}\, 2^{k/2}(k/2)!}, & \text{if $k$ is even.}
   \end{cases}   
\end{align*}
Next, given $\Delta>0$, whose precise value will be chosen later, 
we split $\bXM$ into two disjoint sets
$\bXM^{<}$ and $\bXM^{\ge}$, where
\begin{align*}
  \bXM^{<} &= \Big\{\bx\in\bXM :  
    \Big\vert \frac{m_k(\bx)}{M}-MT(k)\Big\vert <\Delta\Big\}\\
\intertext{ and }
  \bXM^{\ge} &= \Big\{\bx\in\bXM :  
    \Big\vert \frac{m_k(\bx)}{M}-MT(k)\Big\vert \ge\Delta\Big\}.
\end{align*}
Then, the second moment introduced in~\eqref{eqM2odd0} for $k$ odd 
and \eqref{eqM2even} for $k$ even~is
\begin{align*}
    \M(k,M) 
    = \frac{1}{|\bXM|} \sum_{\bx \in \bXM^{<}} \left(\frac{m_k(\bx)}{M} - MT(k)\right)^2
    + \frac{1}{|\bXM|} \sum_{\bx \in \bXM^{\ge}} \left(\frac{m_k(\bx)}{M} - MT(k)\right)^2.
\end{align*}
The first term being positive, it follows
\begin{align*}
    \M(k,M) 
    \ge \frac{1}{|\bXM|} \sum_{\bx \in \bXM^{\ge}} \left(\frac{m_k(\bx)}{M} - MT(k)\right)^2
    \ge \frac{1}{|\bXM|} \sum_{\bx \in \bXM^{\ge}}\Delta^2.
\end{align*}
This implies
\begin{align*}
    \M(k,M) 
    \ge \frac{\Delta^2}{|\bXM|} |\bXM^{\ge}|
    = \frac{\Delta^2}{|\bXM|} 
    \left(|\bXM| - |\bXM^{<}|\right),
\end{align*}
which further implies
\begin{align*}
    \frac{|\bXM^{<}|}{|\bXM|} \ge 1- \Delta^{-2}\M(k,M).
\end{align*}
Using the estimation for the second moment given in Lemma~\ref{LemmaMoment} and choosing $\Delta = M^{-\delta}$, where $\delta>0$ is constant, we obtain:
\begin{align*}
    \frac{1}{|\bXM|}
    \Big\{\bx\in\bXM :  
    \Big\vert \frac{m_k(\bx)}{M}-MT(k)\Big\vert <M^{-\delta}\Big\}
    \ge 1- O\big(M^{2\delta+\eta-1}+ k^2 M^{2\delta-1}\big).
\end{align*}
Note that when $M\to\infty$ while $k$ is fixed,
the above estimate is not trivial if the constants~$\delta$ and~$\eta$ satisfy condition $2\delta+\eta<1$.
Additionally, if the condition holds, $k$ can be allowed to grow 
along with $M$, but not faster 
than~$M^{1/2-\delta}$.
We state the following theorem with the complete result obtained.

%%%%%%%%%%%%%%%%%%%%%%%%%%%%%%%%%%%%%%%
\begin{theorem}\label{Theorem11}
Let $\cM$ be a set of  cardinality $M$, where $M$ is a positive even integer. Let $\psi=\psi_{\cM}$ be an involutive permutation of $\cM$ and consider the set  $\bXM$ of $\psi$-symmetric vectors of length $M$ defined by~\eqref{eqDefinitionbXM}.
Suppose the number of fixed points of $\psi$
is $\theta M\left(1 + O(M^{\eta-1})\right)$, where
$\theta\in [0,1]$ and $\eta\in [0,1)$.
Also suppose $\delta>0$ is a constant that verifies
condition $2\delta+\eta<1$.
Let $k$ be a non-negative integer. 
Then, we have: \\[6pt]
\noindent\textbf{1.} If $k$ is odd, then:
\begin{align*}
   \frac{1}{|\bXM|}
    \Big\{\bx\in\bXM :  
    \Big\vert 
    \frac{m_k(\bx)}{M}- \fdfrac{1-\theta}{ek!}
    \Big\vert <M^{-\delta}\Big\}
    &\ge 1- O\big(M^{2\delta+\eta-1}+ k^2 M^{2\delta-1}\big);  
\end{align*}

\noindent\textbf{2.} 
If $k$ is even, then:
\begin{align*} 
    % &\quad 
    \frac{1}{|\bXM|}
    \Big\{\bx\in\bXM :  
    \Big\vert 
    \frac{m_k(\bx)}{M}- 
     \fdfrac{1-\theta}{e\,k!}
     - &\fdfrac{\theta}{e^{1/2}\, 2^{k/2}(k/2)!}
    \Big\vert <M^{-\delta}\Big\}\\
    &\qquad\quad\,\ge 1- O\big(M^{2\delta+\eta-1}+ k^2 M^{2\delta-1}\big).
\end{align*}
\end{theorem}

%%%%%%%%%%%%%%%%%%%%%%%%%%%%%%%%%%%%%%%
\begin{corollary}\label{Corollary11}
Let us assume that the hypotheses of 
Theorem~\ref{Theorem11} are satisfied.
% In the hypotheses of Theorem~\ref{Theorem1}, when 
Then, the number $m_0(\bx)$ of the elements of $\cM$
that are not represented in $\bx$ is estimated 
for almost all $\bx\in\bXM$ by the next formula:
\begin{align*}
 \frac{1}{|\bXM|}
    \Big\{\bx\in\bXM :  
    \Big\vert 
    \frac{m_0(\bx)}{M}- 
     \fdfrac{1-\theta}{e}
     - \fdfrac{\theta}{e^{1/2}}
    \Big\vert <M^{-\delta}\Big\}
    &\ge 1- O\big(M^{2\delta+\eta-1} \big).
\end{align*}

\end{corollary}

Note that Theorem~\ref{Theorem1} in the introduction is a 
restatement of Theorem~\ref{Theorem11} which is less precise but more intuitive.

%%%%%%%%%%%%%%%%%%%%%%%%%%%%%%%%%%%%%%%%%%%%%%%%%%%%%%
%%%%%%%%%%%%%%%%%%%%%%%%%%%%%%%%%%%%%%%%
% \newpage
\section{Proof of Theorem~\ref{Theorem2} and~\ref{Theorem3}}\label{SectionProofTheorem23}

Here we apply Theorem~\ref{Theorem11} to prove the abundant existence 
of vectors $\bx\in\bXM$ that have have components not only linked by $\psi$
but also belonging to certain sets of quantifiable size.
We remark that while the statement of Theorem~\ref{Theorem2} includes the case 
when $u$ may be missing in $\bx$,
looking at the matter from a slightly different angle, the statement of Theorem~\ref{Theorem3}
provides the threshold size for the existence or absence 
of the highlighted phenomenon.

%%%%%%%%%%%%%%%%%%%%%%%%%%%%%%%%%%%%%%%%%%%%%%%%%
\subsection{Proof of Theorem~\ref{Theorem2}}

Let $K = \max\{r,s\}$,  where we know from the hypothesis 
that~$K \ge 1$. 
Given $\mu > 0$ sufficiently small, 
which will be specified later,
for any integer $k\in\{0,1,\dots,M\}$, 
% with $0\le k\le M$, 
we denote
\begin{equation*}
    \cS(k) := \Big\{
    \bx\in\bXM : 
         \Big\vert m_k(\bx)- \frac{M}{e\,k!}\Big\vert 
         > \frac{\mu}{e\,k!}  M 
    \Big\}
\end{equation*}
and we define the exceptional set  $\cE_K$ 
as the union:
\begin{align*}
    \cE_K := \bigcup_{k=0}^K \cS(k).
\end{align*}
From Theorem~\ref{Theorem11}, under the assumption that 
$\psi$ has no fixed points, it follows that for any~$\varepsilon>0$, there exists an integer \mbox{$M_0 = M_0(k,r,\varepsilon)>0$} 
and $\mu = \mu(\varepsilon)$ that tends to zero when 
$\varepsilon$ tends to zero,
such that for all sets of $\cM$ of cardinality $M \ge M_0$, 
we have
\begin{align}\label{eqSkbound}
    |\cS(k)| \le \frac{\varepsilon}{K+1}|\bXM|,\ \ \text{ for $k=0,1,\dots,M$},
\end{align}
so that it follows
\begin{align}\label{eqEKbound}
    |\cE_K| \le \varepsilon\,|\bXM|.
\end{align}
For any $\bx\in\bXM$ and $k\in[M]$, we let 
\begin{align*}
   \cT_k(\bx) := \{y \in \cM : y \text{ is represented at most $k$ times in } \bx\}.
\end{align*}
Then, for any $\bx\in\bXM\setminus\cE_K$ (which means that
$\bx\not\in\cS(j)$ for any $j=1,2,\dots,k$)
we have
\begin{equation}\label{eqTkbound}
  \begin{split}
   |\cT_k(\bx)| &=
   m_0(\bx)+\cdots+m_k(\bx)\\
   &\ge \Big(\frac{M}{e\,0!}-\frac{\mu}{e\,0!}M\Big)
   +\cdots+
   \Big(\frac{M}{e\,k!}-\frac{\mu}{e\,k!}M\Big)\\
   & = (1-\mu)E(k)M,
  \end{split}  
\end{equation}
where 
\begin{align*}
  E(k) := \frac{1}{e}\sum_{j=0}^k\frac{1}{j!}.  
\end{align*}

Next, given $\mu >0$ and integers $r,s\ge 0$, we let 
$c = c(\varepsilon, r, s)$ be defined~by
\begin{align*}
   c  = \frac{1}{2}\left(3 - (1 - \mu)\big(E(r) + E(s)\big)\right).
\end{align*}
Note that since $\max\{r,s\} \ge 1$, we have:
\begin{align}\label{eqC}
  c \le \frac{1}{2}\left(3 - (1-\mu)\big(E(0) + E(1)\big)\right) \le 
  \frac{1}{2}\left( 3 - \frac{3(1-\mu)}{e}\right) < 1, 
\end{align}
the last inequality taking place if we choose 
$\mu<10^{-2}$, say. 

Now, for any subset $\cB \subset \cM$, which can be large enough such that
$|\cB| > c|\cM|$, on combining~\eqref{eqTkbound} and \eqref{eqC} implies
\begin{equation*}%\label{eqBTbound}
  \begin{split}
  |\cB| + |\cT_k(\bx)| &> c M+ (1-\mu)E(k) M\\ 
   &= \frac{1}{2}\big(3 - (1-\mu)(E(r) + E(s) - 2 E(k)\big) M,
  \end{split}
\end{equation*}
for all $k=0,1,\dots,K$.
 It follows that 
\begin{align}\label{eqBTbound}
    \left|\cB \cap \cT_k(\bx)\right| 
    &\ge  |\cB| + |\cT_k(\bx)| - |\cM|  \notag\\
    &> \frac{1}{2}\big(1 - (1-\mu)(E(r) + E(s) - 2 E(k)\big) M ,
\end{align}
for all $k=0,1,\dots,K$.
And furthermore, for any bijection $\phi:\cM\to\cM$, by applying~\eqref{eqBTbound} 
for $k = r$ and $k = s$, it yields:
\begin{equation*}%\label{eqpsiBTbound}
  \begin{split}
  &\,\quad  \left|\phi\left( B \cap \cT_r(\bx)\right)\right| 
     + \left|\cB \cap \cT_s(\bx)\right| \\
  &= \left|B \cap \cT_r(\bx)\right| 
     + \left|\cB \cap \cT_s(\bx)\right| \\
  &> \frac{1}{2}\big(1 - (1-\mu)(E(s) - E(r)\big) M
   +  \frac{1}{2}\big(1 - (1-\mu)(E(r)- E(s)\big) M\\
  &= M.
  \end{split}
\end{equation*}
As a consequence, we find that
\begin{align*}
   \left|\phi\left( B \cap \cT_r(\bx)\right) 
   \cap \left(\cB \cap  \cT_s(\bx)\right)\right| 
   &\ge  
   \left|\phi\left( B \cap \cT_r(\bx)\right)\right| 
     + \left|\cB \cap \cT_s(\bx)\right| - |\cM| 
     > 0.
\end{align*}
This shows that there exists 
$u \in B \cap \cT_r(\bx)$ such that 
$v = \phi(u) \in \cB \cap \cT_s(\bx)$,  
which completes the proof of Theorem \ref{Theorem2}.
% \end{proof}

%%%%%%%%%%%%%%%%%%%%%%%%%%%%%%%%%%%%%%%%%%%%%%%%%
\subsection{Proof of Theorem~\ref{Theorem3}}

In the case $k=0$ of the unrepresented elements of the set $\cM$
% , of cardinality $M$, 
among the components of~$\bx$
from Theorem \ref{Theorem11}, we know that for each 
$\varepsilon$, say $\varepsilon = 10^{-3}$, 
there exists 
an integer $M_0 = M_0(\varepsilon)$ such that most vectors
$\bx \in \bXM$, precisely, 
at least $(1-\varepsilon)|\bXM|$ of them,  satisfy 
\begin{align}\label{eqm013}
    \left|m_0(\bx) - e^{-1} M\right| < 10^{-3} M. 
\end{align}
% A = numerical_approx(2*(0.87 - exp(-1) - 10^(-3))) #1.00224111765712
% A = numerical_approx(2*(0.86 - exp(-1) - 10^(-3))) #0.982241117657115
% 
% A = numerical_approx(2*(0.87 - exp(-1) - 10^(-2))) #0.984241117657115
% A = numerical_approx(2*(0.86 - exp(-1) - 10^(-2))) #0.964241117657115

Let $\bx = (x_1, \dots, x_M)\in\bXM$ satisfy~\eqref{eqm013}, 
but otherwise arbitrary.
By definition of $\bx$, we see that $M-m_0(\bx)$, which is
the number of $y\in\cM$ that do appear as components of~$\bx$
satisfy $M-Me^{-1}-10^{-3}M\le M-m_0(\bx)$, condition that can be written as
\begin{align}\label{eqT51}
   |\{x_1, \dots, x_M\}| \ge \left(1 - e^{-1} - 10^{-3}\right) M. 
\end{align}
Next, for each sufficiently large subset $\cB\subset\cM$
that has
\begin{align}\label{eqT52}
   |\cB| \ge 0.87 |\cM|
\end{align}
elements, using~\eqref{eqT51}, yields 
\begin{equation}\label{eqT53}
  \begin{split}
 |\cB \cap \{x_1, \dots, x_M\}| &\ge |\cB| + |\{x_1, \dots, x_M\}| - |\cM| \\
 &\ge \left(0.87 - e^{-1} - 10^{-3}\right) M\\
 &>M/2,
   \end{split}
\end{equation}
because $0.87 - e^{-1} - 10^{-3}=0.50112\dots$
% A = numerical_approx((0.87 - exp(-1) - 10^(-3))) #0.501120558828558

Then, an analogous relation can be deduced for any bijection 
$\phi : \cM\to\cM$ and using~\eqref{eqT53}:
\begin{equation}\label{eqT54}
  \begin{split}
    &\quad |\phi(\cB \cap \{x_1, \dots, x_M\}) \cap\left(\cB \cap \{x_1, \dots, x_M\}\right)| \\
    &\ge | \phi(\cB \cap \{x_1, \dots, x_M\})| + |\cB \cap \{x_1, \dots, x_M\}| - |\cM| \\
    &= 2|\cB \cap \{x_1, \dots, x_M\}| - M \\ 
    &\ge  \left( 2 \cdot \left( 0.87 - e^{-1}-10^{-9}\right) - 1\right) M > 0.
  \end{split}
\end{equation}
In other words, the inequality~\eqref{eqT54} implies
the existence of an element $y \in \cB$ such that both~$y$ and $\phi(y)$ are simultaneously among the components of $\bx$.
This completes the proof of Theorem~\ref{Theorem3}.
% \end{proof}

\begin{remark}\label{Remark86Tight}
For any $\bx\in\bXM$ in the statement of Theorem~\ref{Theorem3},
% in the proof of Theorem~\ref{Theorem3}, 
if the subset $\cB\subset\cM$ in the statement of Theorem~\ref{Theorem3} were smaller, containing fewer than~$86\%$ of the elements of $\cM$, then the result would no longer be valid. 
Indeed, in that case, we could take as $\cB$ 
the set of all elements of $\cM$ not represented in $\bx$.
Then we would have
\begin{align*}
   |\cB \cap \{x_1, \dots, x_m\}| \le (0.86  - e^{-1} + 10^{-3})M < M/2,
\end{align*}
% A = numerical_approx((0.86 - exp(-1) + 10^(-3))) #0.493120558828558
and, considering any bijection $\phi : \cM \to \cM$ such that $\phi(\cB \cap \{x_1, \dots, x_M\})$ is entirely contained in $\cM \setminus \left(\cB \cap \{x_1, \dots, x_m\}\right)$,
we would find that there is no 
$y \in \cB \cap \{x_1, \dots, x_m\}$ with $\phi(y) \in \cB \cap \{x_1, \dots, x_m\}$.  

Note that, in the proof of Theorem~\ref{Theorem3},
working with any $\varepsilon$ smaller than~$10^{-3}$,
we see that the exact value of the threshold percentage 
that determines whether the statement 
of Theorem~\ref{Theorem3} remains valid or not is 
$e^{-1}+1/2=0.867879\dots$
% A = numerical_approx((1/2 + exp(-1)))  #0.867879441171442
\end{remark}

%%%%%%%%%%%%%%%%%%%%%%%%%%%%%%%%%%%%%%%%
\section{Closing remarks and further directions}
We have introduced a statistical model to describe the distribution 
of appearances of the elements in a set $\cM$ into finite sequences 
$\bx\in\bXM$, which have length $M=|\cM|$, and whose terms are linked through an involutive bijection $\psi$, as defined in~\eqref{eqDefinitionbXM}.

We have shown that the probability density function of the 
number  $k\ge 0$  that counts the occurrences of components in the sequences coincides with the density in a Poisson process
multiplied by a factor given by the number of fixed points of $\psi$, if $k$ is odd.
To this, an additional different term with a complementary multiplicative factor is added when $k$ is even.
Then we have used the obtained result to show that 
if $\psi$ has no fixed points, then almost all sequences in 
$\bXM$ have terms whose frequency in sequences can be prescribed a priori and, additionally, have complementary interlinks.
While our discussion is based on the assumption that $M$ is presumed to be even, the model can also be extended to the case where $M$ is odd.

To better align with the original intention of applying the model 
to understanding the distribution of the sequence of 
factorials $\pmod p$, 
the set of sequences $\bXM$ can be filtered further, 
retaining in the statistics only those sequences that satisfy the additional non-equal-neighbor hypothesis
(see \cite{CPZ2024} and Broughan and Barnett~\cite{BB2009}).
Additionally, another step in the same direction 
would involve analyzing the model for shorter sequences of length $N$, with $\lambda\sim N/M$ as $M\to\infty$, where $\lambda\in (0,1]$.

%%%%%%%%%%%%%%%%%%%%%%%%%%%%%%%%%%%%%%%%%
\subsection*{Acknowledgments}
We would like to express our thanks to Vladimir Vatutin 
for his prompt and concise presentation 
of the general probabilistic context in which our work fits and for providing 
references~\cite{KSC1978, Mik1980, VM1982}.
%17, 20, 21.
Additionally, we are also grateful to the anonymous referee for such a meticulous review and insightful recommendations for enhancing the presentation of the manuscript.

%%%%%%%%%%%%%%%%%%%%%%%%%%%%%%%%%%%%%%%%%%%%%%
% \bibliographystyle{plainurl}% shows urls
% \bibliography{bibFactorials}

\end{document}